\documentclass{amsart}

\usepackage{hyperref}

\usepackage{amssymb}
\usepackage{amsmath}
\usepackage{amsthm}
\usepackage{enumerate}
\usepackage{graphicx}
\usepackage{tikz-cd}

\theoremstyle{plain}

\newtheorem{theorem}{Theorem}[section]

\newtheorem{proposition}[theorem]{Proposition}
\newtheorem{lemma}[theorem]{Lemma}

\newtheorem{question}[theorem]{Question}
\newtheorem{conjecture}[theorem]{Conjecture}
\newtheorem{problem}[theorem]{Problem}

\theoremstyle{definition}

\newtheorem{example}{Example}[subsection]

\theoremstyle{remark}

\newcommand{\RiemannSphere}{\widehat{\mathbb{C}}}

\DeclareMathOperator{\spann}{span}

\renewcommand{\Re}{\mathrm{Re\,}}

\begin{document}

\title{Analytic Capacity : computation and related problems}

\date{May 7, 2017}

\author[M. Younsi]{Malik Younsi}
\thanks{Supported by NSF Grant DMS-1664807}
\address{Department of Mathematics, University of Washington, Seattle, WA 98195-4350, United States.}
\email{malik.younsi@gmail.com}


\begin{abstract}
We present a brief introduction to analytic capacity, with an emphasis on its numerical computation. We also discuss several related open problems.
\end{abstract}

\maketitle

\section{Introduction}

Analytic capacity is an extremal problem introduced by Ahlfors \cite{AHL} in 1947 in order to study a problem posed by Painlev\'e \cite{PAI} asking for a geometric characterization of the compact plane sets $E \subset \mathbb{C}$ having the property that every bounded analytic function on $\mathbb{C} \setminus E$ is constant. Such compact sets are called \textit{removable} and are precisely the ones of zero analytic capacity. Finding both necessary and sufficient geometric conditions for removability quickly appeared to be very difficult, and it took more than a hundred years until a satisfying solution was obtained, thanks to the work of Melnikov, David, Tolsa and many others.

Another motivation for the study of analytic capacity came to light in 1967, when Vitushkin \cite{VIT} observed that it plays a fundamental role in the theory of uniform rational approximation of analytic functions on compact plane sets. More precisely, analytic capacity can be used to formulate necessary and sufficient conditions for a compact set $E$ to have the property that every function continuous on $E$ and analytic in the interior of $E$ is uniformly approximable by rational functions with poles outside $E$. See e.g. \cite{ZAL} for more details and other applications of analytic capacity to this type of problem.

Although several important problems in Complex Analysis and related areas can be expressed in terms of analytic capacity, the fact that the latter is very hard to calculate or even estimate often remains a major obstacle. Motivated by this and by applications to the subadditivity problem, the author and Ransford developed in \cite{YOU} an efficient and rigorous method based on quadratic minimization for the numerical computation of the analytic capacity of sufficiently nice compact sets. The method yields rigorous upper and lower bounds which converge to the true value of the capacity.

The purpose of this survey article is to present a detailed description of this method and discuss some applications to related open problems. Section \ref{sec2} contains a brief introduction to analytic capacity, including its definition and main properties. Then, in Section \ref{sec3}, we describe the numerical method from \cite{YOU}. The remaining sections are devoted to applications and related open problems, and each of them can be read independently of the others. In Section \ref{sec4}, we study the subadditivity inequality $\gamma(E \cup F) \leq \gamma(E)+\gamma(F)$. More precisely, we use an approximation technique due to Melnikov \cite{MEL} to show that it suffices to prove the inequality in the very special case where $E$ and $F$ are disjoint compact sets which are finite unions of disjoint closed disks, all with the same radius. We also give numerical evidence for subadditivity and formulate a conjecture which, if true, would imply that the inequality holds. Section \ref{sec5} deals with the problem of finding all Ahlfors functions which are rational maps, as instigated by Jeong and Taniguchi \cite{JEO}. Finally, in Section \ref{sec6}, we discuss the relationship between analytic capacity and a similar extremal problem, the Cauchy capacity.

\section{Analytic capacity : definition and preliminaries}
\label{sec2}
In this section, we introduce the notions of analytic capacity, Ahlfors function and Garabedian function. The content is quite standard and can be found in \cite{GAR} for example.

\subsection{Definition and elementary properties}
\label{subsec21}

Let $E$ be a compact subset of the complex plane $\mathbb{C}$ and let $\Omega:= \RiemannSphere \setminus E$ be the complement of $E$ in the Riemann sphere. The \textit{analytic capacity} of $E$ is defined by
$$\gamma(E):= \sup\{ |f'(\infty)|: f \in H^{\infty}(\Omega), \, |f| \leq 1\}.$$
Here $f'(\infty):=\lim_{z \to \infty} z(f(z)-f(\infty))$ denotes the coefficient of $1/z$ in the Laurent expansion of $f$ near infinity, and $H^{\infty}(\Omega)$ denotes the space of all bounded analytic functions on $\Omega$.

The following properties follow more or less directly from the definition.

\begin{itemize}
\item $E \subset F \Rightarrow \gamma(E) \leq \gamma(F)$.
\item $\gamma(aE+b)=|a|\gamma(E) \qquad (a,b \in \mathbb{C})$.
\item Outer regularity : $E_n \downarrow E \Rightarrow \gamma(E_n) \downarrow \gamma(E)$.
\item $E$ is removable for bounded analytic functions if and only if $\gamma(E)=0$.
\item $\gamma(E) = \gamma(\partial_e E)$, where $\partial_e E$ is the outer boundary of $E$, i.e. the boundary of the unbounded component of $\Omega$.

\end{itemize}

Note that the last property implies that the analytic capacity of a set $E$ only depends on its unbounded complementary component. We can therefore assume without loss of generality that $\Omega$ is connected, which we will do for the remaining of this article.

\subsection{The Ahlfors function}
\label{subsec22}

A simple normal families argument shows that for every compact set $E$, there exists an extremal function $f$ for $\gamma(E)$, that is, a function $f$ analytic on $\Omega$ with $|f| \leq 1$ and $f'(\infty)=\gamma(E)$. Moreover, in the case $\gamma(E)>0$, any extremal function $f$ must satisfy $f(\infty)=0$, for otherwise composing with a M\"{o}bius automorphism of the unit disk sending $f(\infty)$ to $0$ would give a larger derivative at $\infty$.

\begin{proposition}
In the case $\gamma(E)>0$, the extremal function is unique.
\end{proposition}
The following proof due to Fisher \cite{FIS} is so short and elegant that we reproduce it below.

\begin{proof}
Suppose that $f$ and $g$ are extremal functions, and let $h:=(f+g)/2$, $k:=(f-g)/2$, so that $f=h+k$ and $g=h-k$. We have to show that $k \equiv 0$. Since $|f|^2 \leq 1$ and $|g|^2 \leq 1$ on $\Omega$, we get that $|h|^2 + |k|^2 \pm 2 \Re h \overline{k} \leq 1$. Adding these two inequalities and dividing by two gives $|h|^2 + |k|^2 \leq 1$, and thus
$$|h| + \frac{1}{2} |k|^2 \leq |h| + \frac{1}{2}(1-|h|^2) \leq |h| + \frac{1}{2} (1+|h|)(1-|h|) \leq |h| + 1-|h| = 1.$$
Assume that $k$ is not identically zero, and write $k^2/2 = a_n/z^n + a_{n+1}/z^{n+1} + \dots$ near $\infty$, where $a_n \neq 0$. Note that $n \geq 2$ since $k(\infty)=0$.

Now, let $\epsilon>0$ be small enough so that $\epsilon |a_n| |z|^{n-1} \leq 1$ on some bounded neighborhood $V$ of $E$, and set $f_1:=h+\epsilon \overline{a_n}z^{n-1}k^2/2$. Then $f_1$ is analytic on $\Omega$ and
$$|f_1| \leq |h| + \epsilon |a_n| |z|^{n-1}\frac{|k|^2}{2} \leq |h| + \frac{1}{2}|k|^2 \leq 1$$
on $V \setminus E$, and thus everywhere on $\Omega$ by the maximum modulus principle. But a simple calculation gives $f_1'(\infty) = h'(\infty) + \epsilon |a_n|^2>\gamma(E)$, a contradiction.

\end{proof}

This unique extremal function is called \textit{the Ahlfors function} for $E$ (or on $\Omega$). It is the unique function $f$ analytic on $\Omega=\RiemannSphere \setminus E$ with $|f| \leq 1$ and $f'(\infty)=\gamma(E)$.

We now mention several well-known properties of the Ahlfors function $f$ for a compact set $E$. We of course assume that $\gamma(E)>0$.

\begin{itemize}
\item Conformal invariance : if $h : \Omega \to \Omega'$ is conformal with expansion $F(z)=a_1z+a_0+a_{-1}/z+\dots$ near $\infty$, and if $E':=\RiemannSphere \setminus \Omega'$, then $g:=f \circ h^{-1}$ is the Ahlfors function for $E'$ and $\gamma(E')=|a_1|\gamma(E)$.
\item If $E$ is connected, then $f : \Omega \to \mathbb{D}$ is the unique conformal map of $\Omega$ onto the open unit disk $\mathbb{D}$ normalized by $f(\infty)=0$ and $f'(\infty)>0$. In particular, we have $\gamma(\overline{\mathbb{D}}(z_0,r))=r$ and $\gamma([a,b])=(b-a)/4$.
\item More generally, if $\Omega$ is a non-degenerate $n$-connected domain (i.e. $E$ has $n$ connected components, each of them containing more than one point), then $f:\Omega \to \mathbb{D}$ is a degree $n$ proper analytic map.
\item If $E \subset \mathbb{R}$, then the Ahlfors function for $E$ is given by
$$f(z) = \frac{e^{h(z)}-1}{e^{h(z)}+1} \qquad (z \in \Omega),$$
where
$$h(z):= \frac{1}{2} \int_E \frac{dt}{z-t}.$$
\end{itemize}

The first two properties follow directly from a simple change of variable and Schwarz's lemma respectively. Note that the second property implies that the analytic capacity of a connected set is equal to its logarithmic capacity. The third property was proved by Ahlfors \cite{AHL}. Finally, the last property follows from a result of Pommerenke (see e.g. \cite[Chapter 1, Theorem 6.2]{GAR}) saying that the capacity of a set $E \subset \mathbb{R}$ equals a quarter of its one-dimensional Lebesgue measure $m(E)$. Indeed, one can easily check that the given function $f$ maps $\Omega$ into $\mathbb{D}$ and satisfies $f'(\infty)=m(E)/4$, so it has to be the Ahlfors function.

\subsection{The Garabedian function}
\label{subsec23}

As observed by Garabedian \cite{GARA}, the analytic capacity of sufficiently nice compact sets can be obtained as the solution to a dual extremal problem.

More precisely, let $E \subset \mathbb{C}$ be a compact set, and assume that $E$ is bounded by finitely many mutually exterior analytic Jordan curves. Let $A(\Omega)$ denote the space of analytic functions in $\Omega=\RiemannSphere \setminus E$ which extend continuously up to $\overline{\Omega}$. The \textit{Garabedian function} for $E$ (or on $\Omega$) is the unique function $\psi \in A(\Omega)$ satisfying $\psi(\infty)=1/2\pi i$ and
$$\int_{\partial E} |\psi(\zeta)| \, |d\zeta| = \inf \left\{ \int_{\partial E} |h(\zeta)| \, |d\zeta| : h \in A(\Omega), \, h(\infty)=\frac{1}{2\pi i} \right\}.$$
The Garabedian function has the following properties.

\begin{itemize}
\item $\psi$ extends analytically across $\partial \Omega$.
\item $\psi$ represents evaluation of the derivative at $\infty$, in the sense that, for all $g \in A(\Omega)$, we have
$$g'(\infty) = \int_{\partial E} g(\zeta) \psi(\zeta) \, d\zeta.$$
\item $\psi$ has a well-defined analytic logarithm in $\Omega$. In particular, there exists a function $q \in A(\Omega)$ such that $q(\infty)=1$ and $q^2=2\pi i \psi$.
\item $\gamma(E)=\int_{\partial E} |\psi(\zeta)| \, |d\zeta|.$
\item The differential $f(\zeta)\psi(\zeta) \, d\zeta$ is positive on $\partial E$.
\end{itemize}

We also mention that the above function $q$ is, up to a multiplicative constant, the Szeg\"{o} kernel function $K(\cdot,\infty)$ (the reproducing kernel for the evaluation functional at $\infty$ in the Hardy space $H^2(\Omega)$). Indeed, we have

\begin{equation}
\label{eqszego}
K(z,\infty) = \frac{1}{2\pi \gamma(E)} q(z) \qquad (z \in \Omega).
\end{equation}
See \cite[Theorem 4.3]{GAR}.

\section{Numerical computation}
\label{sec3}

In this section, we present the method developed in \cite{YOU} for the numerical computation of analytic capacity.

First, for the sake of completeness, let us briefly describe a simple method for the computation of $\gamma(E)$ based on Equation (\ref{eqszego}).

It is well-known that if $E$ is bounded by finitely many mutually exterior analytic Jordan curves, then the rational functions with poles outside $E$ are dense in the Hardy space $H^2(\Omega)$. Applying the Gram-Schmidt procedure, one can therefore extract an orthonormal basis $\{u_n\}_{n \geq 1}$. Then, since
$$K(z,\infty) = \sum_{n \geq 1} \overline{u_n}(\infty) u_n(z),$$
we obtain, from Equation (\ref{eqszego}),

$$\gamma(E) = \frac{1}{2\pi} K(\infty,\infty)^{-1} = \frac{1}{2\pi} \left( \sum_{n \geq 1} |u_n(\infty)|^2 \right)^{-1}.$$

However, this very simple method has two main disadvantages. First, computing a large number of elements of the orthonormal basis $\{u_n\}_{n \geq 1}$ seems to be numerically unstable, which causes precision issues. Secondly, and more importantly, truncating the sum in the above formula only yields an upper bound for $\gamma(E)$.

We also mention that the Szeg\"{o} kernel (and therefore the analytic capacity) can be computed by solving the Kerzman-Stein integral equation, see \cite{BEL} and \cite{BOL}. However, as far as we know, there seems to be very few numerical examples, and only for compact sets with a few number of connected components.

On the other hand, the method developed in \cite{YOU} for the numerical computation of analytic capacity is both efficient and precise for a wide variety of compact sets. Furthermore, it yields both lower and upper bounds which can be made arbitrarily close to the true value of the capacity, thereby providing adequate error control. The method relies on two key estimates, which we present in the next subsection.

\subsection{Estimates for analytic capacity}
\label{subsec31}

Let $E \subset \mathbb{C}$ be a compact set, and again suppose that $E$ is bounded by finitely many mutually exterior analytic Jordan curves.

\begin{theorem}[Younsi--Ransford \cite{YOU}]
\label{thmest}
We have
\begin{equation}
\label{estime1}
\gamma(E) = \min \left\{ \frac{1}{2\pi}\int_{\partial E} |1+g(z)|^2|dz|: g \in A(\Omega), \, g(\infty)=0 \right\}
\end{equation}
and
\begin{equation}
\label{estime2}
\gamma(E) = \max \left\{ 2\Re h'(\infty) - \frac{1}{2\pi}\int_{\partial E} |h(z)|^2 |dz|: h \in A(\Omega), \, h(\infty)=0 \right\}.
\end{equation}
Here the minimum and maximum are attained respectively by the functions $g=q-1$ and $h=fq$, where $f$ is the Ahlfors function for $E$ and $q$ is the square root of $2\pi i$ times the Garabedian function $\psi$ for $E$, as in Subsection \ref{subsec23}.

\end{theorem}

The first estimate is due to Garabedian and follows from the formula for $\gamma(E)$ in terms of $\psi$, as in Susbsection \ref{subsec23}. The second estimate was obtained in \cite{YOU}. We reproduce its proof for the reader's convenience.

\begin{proof}
A simple calculation shows that the maximum in (\ref{estime2}) is attained by the function $fq$, which belongs to $A(\Omega)$ and vanish at $\infty$. Note that $|f|=1$ on $\partial E$. It thus suffices to show that if $h \in A(\Omega)$ and if $h(\infty)=0$, then

$$ \gamma(E) \geq 2\Re h'(\infty) - \frac{1}{2\pi}\int_{\partial E} |h(z)|^2 |dz|.$$
Denote by $T(z)$ the unit tangent vector to $\partial E$ at $z$,
that is, $dz = T(z) |dz|$ with $|T| \equiv 1$.
Let $\langle h_1, h_2 \rangle$ denote
$\int_{\partial E} h_1(z) \overline{h_2(z)} |dz|$
and $\|h\|_2^2:=\langle h,h \rangle$.
Then
$$
0 \leq \frac{1}{2\pi}\|h-i\overline{qT}\|_{2}^2
= \frac{1}{2\pi} \|h\|_{2}^2 + \frac{1}{2\pi} \|q\|_{2}^2 + 2\Re \frac{1}{2\pi} \langle h,-i\overline{qT} \rangle,
$$
so that
$$
0 \leq  \frac{1}{2\pi}\int_{\partial E} |h(z)|^2|dz|
+ \gamma(E)  -2 \Re  \frac{1}{2\pi i}\int_{\partial E} h(z)q(z)dz.
$$
It follows that
$$
0 \leq  \frac{1}{2\pi}\int_{\partial E} |h(z)|^2|dz|
+  \gamma(E)  -2 \Re (hq)'(\infty).
$$
Since $q(\infty)=1$ and $h(\infty)=0$, we have $(hq)'(\infty) = h'(\infty)$, and thus
$$
\gamma(E) \geq   2 \Re h'(\infty) - \frac{1}{2\pi}\int_{\partial E} |h(z)|^2 |dz|,
$$
as required.

\end{proof}

\subsection{The algorithm}
\label{subsec32}

Theorem \ref{thmest} yields a simple method for the numerical computation of analytic capacity based on quadratic minimization. We only explain how to obtain a decreasing sequence of upper bounds for $\gamma(E)$, the method for the lower bounds being very similar.

Let $S$ be any finite set containing at least one point in each component of the interior of $E$. It is well-known that the set of rational functions with poles in $S$ is uniformly dense in $A(\Omega)$, by Mergelyan's theorem for example.

For each $k \in \mathbb{N}$, let $\mathcal{F}_k$ be the set of all functions of the form $(z-p)^{-j}$, with $1 \leq j \leq k$ and $p\in S$. In view of (\ref{estime1}), the quantity
$$
u_k := \min \left\{\, \frac{1}{2\pi}\int_{\partial E} |1+g(z)|^2|dz| : g \in \spann \mathcal{F}_k \,\right\}
$$
gives an upper bound for $\gamma(E)$.

In practice, we find the minimum $u_k$ as follows.

\begin{enumerate}
\item We define
$$
g(z):= \sum_{g_j \in \mathcal{F}_k} (c_j + id_j) g_j(z),
$$
where the $c_j$'s and $d_j$'s are real numbers to determine.

\item We compute the integral
$$
\frac{1}{2\pi}\int_{\partial E} |1+g(z)|^2|dz|
$$
as an expression in the $c_j$'s and $d_j$'s. This gives a quadratic form with a linear term and a constant term.

\item We find the $c_j$'s and $d_j$'s that minimize this expression, by solving the corresponding linear system.
\end{enumerate}

The sequence $(u_k)_{k=1}^\infty$ is decreasing by construction, and it converges to $\gamma(E)$ by Mergelyan's theorem.

\newpage

\subsection{Numerical examples}
\label{subsec33}

In this subsection, we present several numerical examples to illustrate the method. All the numerical work was done with \textsc{matlab}.

\begin{example}
\label{ex1}

{\em Union of two disks.}

Here $E$ is the union of two disks of radius $1$ centered at $-2$ and $2$.

\begin{figure}[h!t!b]
\begin{center}
\includegraphics[width=7cm, height=7cm]{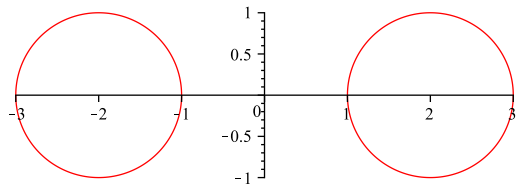}
\caption{The compact set $E$ for Example \ref{ex1}}
\end{center}
\end{figure}

In this case, a natural choice for the set $S$ described in Subsection \ref{subsec32} is $\{-2,2\}$, so that lower and upper bounds for $\gamma(E)$ can be computed using linear combinations of functions of the form $(z \pm 2)^{-j}$. However, we shall instead consider functions of the form
$$
g_j(z)=\frac{1}{z-a_j},
$$
where the $a_j$'s are distinct points in the interior of $E$. The reason is that with simples poles instead of multiple ones, the integrals involved in the method are easily calculated analytically, which results in a significant gain in efficiency.

Typically, for each disk centered at $c$ with radius $r$, we put the poles $a_j$ at the points
$$
\{c, \,c\pm r_1, \,c\pm r_1i, \,c\pm r_2, \,c\pm r_2i, \dots, c\pm r_n,\,c\pm r_ni\},
$$
where $r_1, \dots, r_n$ are equally distributed between $0$ and $r$. This choice, although completely arbitrary, seems to yield good numerical results, as shown in Table~\ref{table0} below.
\newpage

\begin{table}[!hbp]
\begin{center}
\caption{Lower and upper bounds for $\gamma(E)$ for Example \ref{ex1}}
\label{table0}
\begin{tabular}{|c|l|l|r|}
\hline
Poles per disk & Lower bound for $\gamma(E)$ & Upper bound for $\gamma(E)$ & Time (s) \\
\hline
$1$ & 1.875000000000000 &  1.882812500000000 & 0.001867 \\
$5$ & 1.875593064023693 &  1.875619764386366 & 0.003115 \\
$9$ &  1.875595017927203 &  1.875595038756883 & 0.003462 \\
$13$ & 1.875595019096871 &  1.875595019097141 &  0.003854 \\
$17$ & 1.875595019097112 & 1.875595019097164 & 0.005046 \\
\hline
\end{tabular}
\end{center}
\end{table}

The numerical results agree with the value
$$
\gamma(E)\approx 1.8755950190971197289,
$$
obtained using the formula from \cite{MUR2} for the analytic capacity of the union of two disks in terms of elliptic integrals.
\end{example}

\begin{example}
\label{ex2}
{\em Union of several disks.}

In this example, we compute the analytic capacity of the union of several closed disjoint disks, all with the same radius. Such compact sets are especially relevant for the subadditivity problem, as we will see in Section \ref{sec4}.

The centers of the disks as well as the radius were randomly generated, and we computed the analytic capacity using $5$ poles per disk.

\begin{table}[!hbp]
\begin{center}
\caption{Lower and upper bounds for $\gamma(E)$ for Example \ref{ex2}}
\label{table1}
\begin{tabular}{|c|l|l|r|}
\hline
Number of disks & Lower bound for $\gamma(E)$ & Upper bound for $\gamma(E)$ & Time (s) \\
\hline
$10$ & 4.565899720026281 &  4.565918818251154 & 0.008976 \\
$100$ & 0.856133344195575 &  0.856133345307269 & 1.465479 \\
$200$ & 4.968896699483210 & 4.968903654028726 & 10.525317 \\
$500$ & 6.471536902032636 & 6.471540680554522 & 163.594474 \\
$1000$ & 6.548856893117339 & 6.548862005607368 & 1252.940712\\
\hline
\end{tabular}
\end{center}
\end{table}
\end{example}
\newpage
\begin{example}
\label{ex30}
{\em Union of four ellipses.}

The compact set $E$ is composed of four ellipses centered at $-3$, $3$, $10i$, $-10i$. Each ellipse has a semi-major axis of $2$ and a semi-minor axis of $1$.

\begin{figure}[!h]
\begin{center}
\includegraphics[width=8cm, height=8cm]{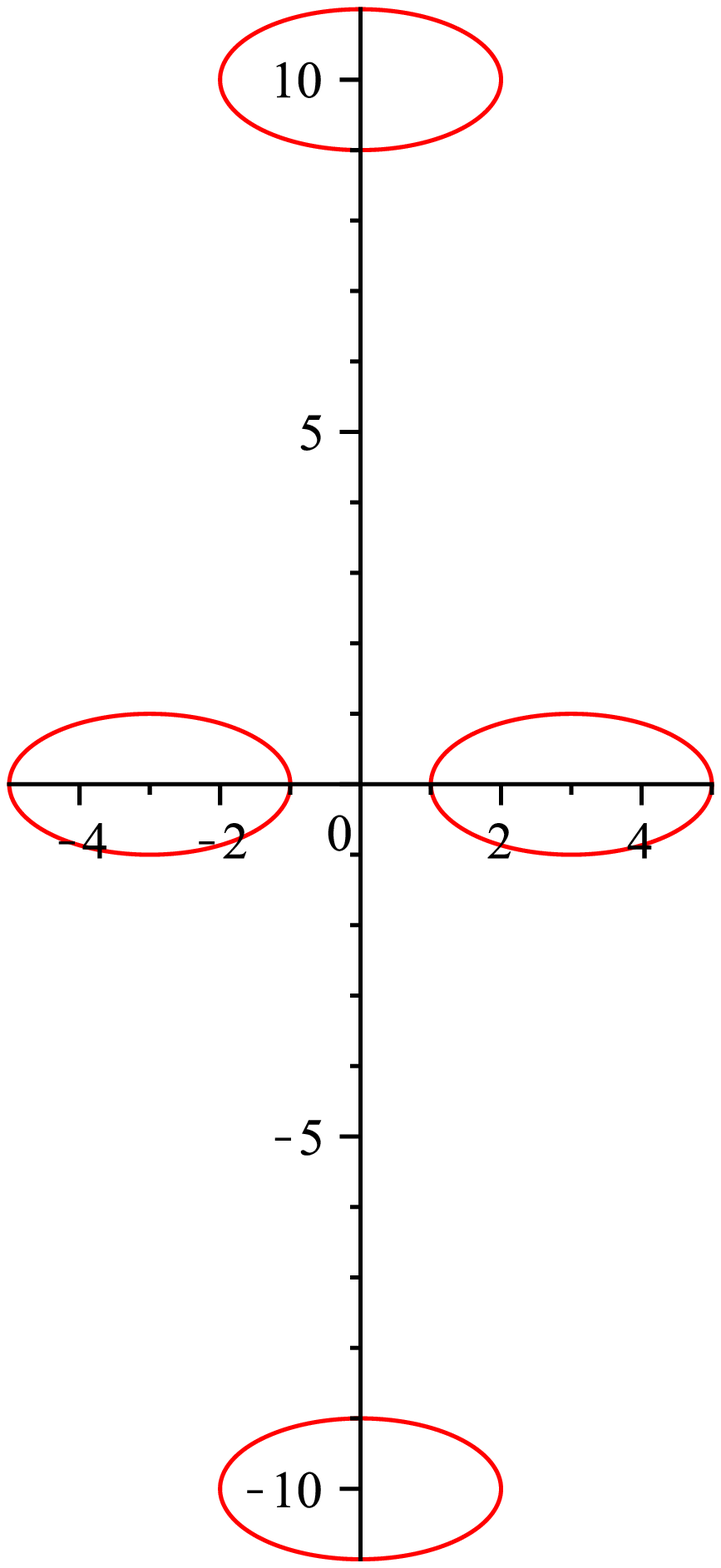}
\caption{The compact set $E$ for Example \ref{ex30}}
\end{center}
\end{figure}

\begin{table}[!h]
\caption{Lower and upper bounds for $\gamma(E)$ for Example \ref{ex30}}
\begin{center}
\begin{tabular}{|c|l|l|r|}
\hline
Poles per ellipse & Lower bound for $\gamma(E)$ & Upper bound for $\gamma(E)$ & Time (s) \\
\hline
$1$ &   4.290494449193028 &      5.652385361295098&     0.962078 \\
$5$ &   5.252560204660928 &      5.409346641724527 &    17.268477 \\
$9$ &   5.356419530523225 &     5.377445892435984 &    54.260216 \\
$13$ &   5.370292494009306 &    5.372648058950175 &  111.424592 \\
$17$ &    5.371877137036634 &     5.372044462730262&   190.042871 \\
$41$ &    5.371995432221965 &   5.371995878776166 & 1100.468881 \\
\hline
\end{tabular}
\end{center}
\end{table}

In this case, the integrals involved have to be calculated numerically. We used a recursive adaptive Simpson quadrature with an absolute error tolerance of $10^{-9}$.
\end{example}

\subsection{Further numerical examples : piecewise analytic boundary}
\label{subsec34}

It was proved in \cite{YOU} that the estimates of Theorem \ref{thmest} remain valid in the case of compact sets bounded by finitely many piecewise-analytic curves, provided the space $A(\Omega)$ is replaced by a larger space of analytic functions, the Smirnov space $E^2(\Omega)$. The convergence of the lower and upper bounds to the capacity was also established in this more general case.

As a consequence, given a compact set $E$ bounded by finitely many mutually exterior piecewise-analytic curves, the algorithm of Section \ref{subsec32} should in principle give an approximation of $\gamma(E)$.

\begin{example}
\label{ex3}
{\em The square.}

In this example, we consider the square with corners $1$, $i$, $-1$, $-i$.

\begin{figure}[!h]
\begin{center}
\includegraphics[width=5cm, height=5cm]{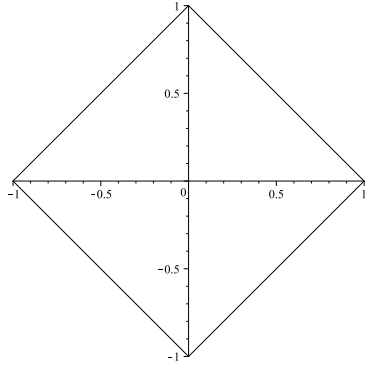}
\caption{The compact set $E$ for Example \ref{ex3}}
\end{center}
\end{figure}

We computed the bounds using functions of the form $z^{-j}$ for $1 \leq j \leq k$. Table~\ref{Ta:square} lists the bounds obtained for different values of $k$.

\begin{table}[!h]
\begin{center}
\caption{Lower and upper bounds for $\gamma(E)$ for Example \ref{ex3}}
\label{Ta:square}
\begin{tabular}{|c|l|l|r|}
\hline
$k$ & Lower bound for $\gamma(E)$ & Upper bound for $\gamma(E)$ & Time (s) \\
\hline
$2$ &  0.707106781186547 &     0.900316316157106 &     0.021981 \\
$3$ &   0.707106781186547 &   0.900316316157106 & 0.069278\\
$4$ &  0.707106781186547 &   0.887142803070031 &   0.109346\\
$5$ & 0.746499705182962 &   0.887142803070031 &    0.145614 \\
$6$ &   0.746499705182962 & 0.887142803070031 &  0.202309\\
$7$ &   0.746499705182962 &     0.887142803070031 &     0.295450 \\
$8$ &   0.746499705182962 &  0.881014562149127 & 0.347996\\
$9$ &   0.761941423753061 &    0.881014562149127 &    0.414684\\
$10$ &  0.761941423753061 &  0.881014562149127 &     0.595552 \\
$15$ &  0.770723484232218 &  0.877175902241141 &      2.425285 \\
$20$ &  0.776589045256849 &  0.872341829081944 &     5.537981 \\
$25$ &  0.784189460107018 &   0.870656623669828 &      10.002786 \\
$30$ &  0.786857803378602 &   0.869257904380382 &      16.344379 \\
$35$ &  0.789068961951613 &    0.868068649269412 &      26.109797 \\
$40$ &   0.790942498354322 &    0.866133165258689 &      33.595790 \\
\hline
\end{tabular}
\end{center}
\end{table}
\end{example}

We notice that the convergence is very slow, especially compared to the results obtained in the analytic boundary case. The main issue here is that the approximation process does not take into account the geometric nature of the boundary. However, as observed in \cite{YOU}, this is easily fixed as follows.

Recall from Theorem \ref{thmest} that the functions to be approximated are $q$ and $fq$, where $f$ is the Ahlfors function for $E$ and $q$ is the square root of $2\pi i$ times the Garabedian function. One can show that these functions behave like $\sqrt{F'}$ near the corners of the square $E$, where $F$ is a conformal map of $\Omega=\RiemannSphere \setminus E$ onto $\RiemannSphere \setminus \overline{\mathbb{D}}$. But if $a$ is one of those corners, then $F$ behaves like $(z-a)^{2/3}$ near $a$, and thus $\sqrt{F'}$ behaves like $(z-a)^{-1/6}$. Since we want functions that are analytic near $\infty$, we consider instead $(1-a/z)^{-1/6}$. This shows that using linear combinations of the functions
$$
\frac{f_j(z)}{z^k}
$$
for $j=0,1,2,3,4$ and $k=1,2, \dots, n$, where $a_1,a_2,a_3,a_4$ are the corners of the square,
$$
f_0(z):=1
$$
and
$$
f_j(z):= \Bigl( \frac{z-a_j}{z} \Bigr)^{-1/6} \qquad (j=1,2,3,4),
$$
should yield better results.

We use this improvement of the method to recompute the analytic capacity of the square. One can see that the convergence is indeed significantly faster.

\begin{table}[!h]
\label{tab1}
\begin{center}
\caption{Lower and upper bounds for $\gamma(E)$}
\begin{tabular}{|c|l|l|r|}
\hline
$n$  & Lower bound for $\gamma(E)$ & Upper bound for $\gamma(E)$ & Time (s) \\
\hline
$2$  & 0.834566926465074 &    0.835066810881929 &     1.334885 \\
$3$  &   0.834609482283050 & 0.834678782816948 & 2.918624\\
$4$  &  0.834622127643984 &  0.834628966618492 &  5.220941\\
$5$  & 0.834626255962448 &  0.834627566559480 &   8.022274 \\
$6$ &  0.834626584020641 & 0.834627152182154 &  11.542859\\
\hline
\end{tabular}
\end{center}
\end{table}

In this case, the answer can be calculated exactly. Indeed, since $E$ is connected, we have that
$$
\gamma(E)
= \operatorname{cap}(E)
= \sqrt{2} \frac{\Gamma(1/4)^2}{4\pi^{3/2}} \approx 0.83462684167407318630,
$$
where $\operatorname{cap}(E)$ is the logarithmic capacity of $E$.
\newpage
Lastly, the improved method can easily be generalized to compute the analytic capacity of any compact set with piecewise analytic boundary. Here is a non-polygonal example.

\begin{example}
\label{ex4}
{\em Union of a disk and two semi-disks}

\begin{figure}[!h]
\begin{center}
\includegraphics[width=7cm, height=7cm]{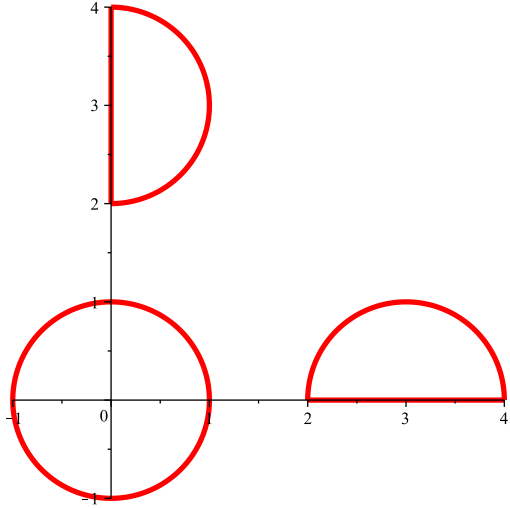}
\caption{The compact set $E$ for Example \ref{ex4}}
\end{center}
\end{figure}

The compact set $E$ is composed of the unit disk and two half-unit-disks centered at $3$ and $3i$.

\begin{table}[h]
\begin{center}
\caption{Lower and upper bounds for $\gamma(E)$ for Example \ref{ex4}}
\begin{tabular}{|c|l|l|r|}
\hline
$n$ & Lower bound for $\gamma(E)$ & Upper bound for $\gamma(E)$ & Time (s) \\
\hline
$2$ & 2.118603690751346 &    2.123888275897654 &     2.546965 \\
$3$  &    2.120521869940459 & 2.121230615594293 &  4.926440\\
$4$  &   2.120666182274863 &  2.120803766391281 &  9.488024 \\
$5$  &  2.120694837101383 &   2.120716977856280 &   13.679742 \\
$6$  &   2.120703235395670 &  2.120709388805280&  22.344576\\
$7$  &   2.120704581010457 & 2.120707633546616 &  28.953791\\
$8$  &   2.120705081159854 &  2.120706704970516&  34.781046\\
\hline
\end{tabular}
\end{center}
\end{table}

\end{example}

\section{The subadditivity problem}
\label{sec4}

This section deals with the subadditivity problem for analytic capacity.

The study of the analytic capacity of unions of sets was instigated by Vitushkin \cite{VIT}, who conjectured that analytic capacity is \textit{semi-additive}, i.e.
$$\gamma(E \cup F) \leq C(\gamma(E) + \Gamma(F)) \qquad (E,F \subset \mathbb{C} \, \mbox{compact} \, ),$$
motivated by applications to problems of uniform rational approximation of analytic functions.
The semi-additivity of analytic capacity was finally proved many years later by Tolsa \cite{TOL3}, along with a solution to Painlev\'e's problem. See \cite{TOL}, \cite{TOL2} and also Subsection \ref{subsec61} for more information regarding Tolsa's work on analytic capacity.

Despite such remarkable work, the optimal value of the constant $C$ remains unknown. In particular, can we take $C=1$?

\begin{problem}
\label{prob1}
Is analytic capacity subadditive? In other words, is it true that
\begin{equation}
\label{subaeq}
\gamma(E \cup F) \leq \gamma(E)+\gamma(F)
\end{equation}
for all compact sets $E,F$?
\end{problem}

We mention that (\ref{subaeq}) is known to hold in the following cases :

\begin{itemize}
\item The sets $E,F$ are disjoint and connected \cite{SUI}.
\item One of the two sets have zero analytic capacity, see e.g. \cite[Lemma 2.6]{YOU2}.
\item The sets $E,F$ are contained in $\mathbb{R}$, since in this case analytic capacity is equal to a quarter of Lebesgue measure, by Pommerenke's result (see the last property in Subsection \ref{subsec22}).
\end{itemize}

Now, we saw in Subsection \ref{subsec33} that the numerical method of \cite{YOU} for the computation of analytic capacity is very efficient in the case of compact sets bounded by analytic curves. Therefore, from a computational point of view, it would certainly be interesting to know whether it suffices to prove (\ref{subaeq}) only for such sets. Using a discrete approach to analytic capacity due to Melnikov \cite{MEL}, one can actually prove a much stronger result.

\begin{theorem}
\label{disks}
The subadditivity of analytic capacity is equivalent to
$$\gamma(E \cup F) \leq \gamma(E)+\gamma(F)$$
for all disjoint compact sets $E$ and $F$ that are finite unions of disjoint closed disks, all with the same radius.
\end{theorem}

\subsection{Melnikov's approach}
\label{subsec41}

In this subsection, we give a proof of Theorem \ref{disks} using an approximation technique due to Melnikov. The proof can also be found in \cite[Theorem 7.1]{YOU}.

Following \cite{MEL}, we introduce the following notation. Let $z_1, \dots, z_n \in \mathbb{C}$ and let $r_1,\dots,r_n$ be positive real numbers. Set $Z:=(z_1,\dots,z_n)$ and $R:=(r_1,\dots,r_n)$. We assume that $|z_j-z_k| > r_j + r_k$ for $j \neq k$, so that the closed disks $\overline{\mathbb{D}}(z_j,r_j)$ are pairwise disjoint. Define $K(Z,R):=\cup_{j=1}^n \overline{\mathbb{D}}(z_j,r)$ and $\mu_1(Z,R):=\sup \{ |\sum_{j=1}^n a_j|\}$, the supremum being taken over all points $a_1, \dots, a_n \in \mathbb{C}$ such that
$$\left| \sum_{j=1}^n \frac{a_j}{z-z_j} \right| \leq 1 \qquad (z \notin K(Z,R)).$$
Clearly, we have $\mu_1(Z,R) \leq \gamma(K(Z,R))$.

Finally, for any compact set $K \subset \mathbb{C}$ and $\delta>0$, we denote the closed $\delta$-neighborhood of $K$ by $K_\delta$.

The proof of Theorem \ref{disks} relies on the following lemma due to Melnikov \cite[Lemma 1]{MEL}.

\begin{lemma}[Melnikov \cite{MEL}]
\label{LemMelnikov}
Let $K \subset \mathbb{C}$ be compact, and let $\delta, \epsilon >0$. Then there exist $z_1,\dots,z_n \in K_{\delta}$ and $0<r<\delta$ such that $|z_j-z_k|>2r$ for $j \neq k$, and
$$\mu_1(Z,R) \geq (1-\epsilon)\gamma(K),$$
where $Z=(z_1,\dots,z_n)$ and $R=(r,\dots,r)$.

In particular,
$$\gamma(K(Z,R)) \geq (1-\epsilon)\gamma(K).$$
\end{lemma}
In other words, the analytic capacity of any compact set can be approximated by the capacity of the union of finitely many disjoint closed disks, all with the same radius (note that if $\delta>0$ is sufficiently small, then $\gamma(K_\delta)$ is close to $\gamma(K)$ by outer-regularity).

We now give a sketch of the proof of Lemma \ref{LemMelnikov}.

\begin{proof}
The idea is to choose finitely many complex numbers $a_1, \dots, a_n$ such that $\gamma(K)=\sum_{j=1}^n a_j$ and
$$\left|\sum_{j=1}^n \frac{a_j}{z-z_j} \right| \leq 1+ C\epsilon \qquad (z \notin K(Z,R)),$$
where $K(Z,R)$ is the union of $n$ disjoint closed disks near $K$. Clearly, this implies the conclusion of the lemma.

In the following, the letter $C$ will denote a positive constant independent of $\delta$ and $\epsilon$, whose value may change throughout the proof.

Let $\phi$ be a $C^\infty$ function on $\mathbb{C}$ with $0 \leq \phi \leq 1$, $\phi = 0$ on $K_{\delta/3}$, $\phi =1$ on $\mathbb{C} \setminus K_{(2/3)\delta}$ and $|\overline{\partial} \phi| \leq C \delta^{-1}$. If $f$ is the Ahlfors function for $K$, then a simple application of the Cauchy-Green formula gives
\begin{equation}
\label{eqCG}
f(z)\phi(z) = \frac{1}{\pi} \int \frac{f(w) \overline{\partial}\phi(w)}{z-w} \, dA(w) \qquad (z \notin K)
\end{equation}

Now, partition the plane with squares $\{Q_j\}$ of side-length $\rho$, where $\rho>0$ is sufficiently small and will be determined later. Note that if $\rho < \delta/6$ and $Q_j \cap \operatorname{supp}(\overline{\partial}\phi) \neq \emptyset$, then $Q_j$ is contained in $K_\delta \setminus K$. Let $J$ be the set of all $j$ for which $Q_j \cap \operatorname{supp}(\overline{\partial}\phi) \neq \emptyset$, so that $J$ is finite. For $j \in J$, define
$$f_j(z):= \frac{1}{\pi} \int_{Q_j} \frac{f(w) \overline{\partial}\phi}{z-w} \, dA(w).$$
Then by Equation (\ref{eqCG}), we have $f \phi = \sum_{j \in J} f_j$ on $\mathbb{C} \setminus K$.
Also, let
$$a_j := \frac{1}{\pi} \int_{Q_j} f(w) \overline{\partial}\phi \, dA(w).$$
Then $|a_j| \leq C \delta^{-1} \rho^2$. Also, we have, for $j \in J$ and $z \notin K$,
$$|f_j(z)| \leq C \delta^{-1} \rho,$$
by splitting the integral over $Q_j \cap \{w:|z-w| \leq \rho\}$ and $Q_j \cap \{w:|z-w|>\rho\}$.

Now, note that if $z_j$ is the center of $Q_j$ and if $|z-z_j| > \rho$, then
\begin{eqnarray*}
\left| f_j(z) - \frac{a_j}{z-z_j} \right| &\leq& \frac{1}{\pi} \int_{Q_j} \frac{|f(w) \overline{\partial}\phi(w)||z_j-w|}{|z-w||z-z_j|} \, dA(w) \\
&\leq& C \delta^{-1} \frac{1}{|z-z_j|} \int_{Q_j} \frac{|z_j-w|}{|z-w|} \, dA(w).
\end{eqnarray*}

To estimate this last integral, note that if $w \in Q_j$, then
$$|z-z_j||w-z_j| \leq C\rho|z-w|,$$
which gives
$$\left| f_j(z) - \frac{a_j}{z-z_j} \right| \leq C \delta^{-1} \rho^3 |z-z_j|^{-2}.$$
Let $r:=N \rho^2 \delta^{-1} \epsilon^{-1}$, where $N$ is a large positive integer to be selected below. We shall take $\rho$ sufficiently small so that $\overline{\mathbb{D}}(z_j,r) \subset Q_j \subset K_\delta \setminus K$ for each $j$. If $z \in \cup_{j \in J} \partial \mathbb{D}(z_j,r)$, then we have
\begin{eqnarray*}
\left| \sum_{j \in J} \frac{a_j}{z-z_j} \right| &\leq& \left| \sum_{j \in J} \left( f_j(z)-\frac{a_j}{z-z_j} \right) \right| + |f(z)\phi(z)|\\
&\leq& 1+\sum_{j \in J'} (C \delta^{-1}\rho + C \delta^{-1}r^{-1}\rho^2) + \sum_{j \in J''} C \delta^{-1} \rho^3 |z-z_j|^{-2},
\end{eqnarray*}
where $J'$ denotes all the $j \in J$ such that $|z-z_j| \leq \rho$, and $J''$ the remaining indices in $J$. We have to show that each of the two sums is small, say less than $\epsilon$. First, note that $J'$ contains at most four indices, and thus the sum over $j \in J'$ can be made smaller than $\epsilon$ by choosing $\rho$ sufficiently small, provided
$$C N^{-1}\epsilon = C \delta^{-1}r^{-1}\rho^2$$
is sufficiently small, which will be the case for $N$ much larger than $C$.

We now estimate the sum over $j \in J''$. Note that if $j \in J''$ and $w \in Q_j$, then
$$|z-w| \leq |z-z_j|+|z_j-w| \leq (1+\sqrt{2}/2)|z-z_j|,$$
so that
$$\sum_{j \in J''} C \delta^{-1} \rho^3 |z-z_j|^{-2} \leq C \delta^{-1} \rho \sum_{j \in J''} \int_{Q_j} \frac{dA(w)}{|w-z|^2}.$$
Now, again if $w \in \cup_{j \in J''}Q_j$, then $|w-z| \geq (1-\sqrt{2}/2)\rho$ and $|w-z| \leq \operatorname{diam}(K_\delta)$. It follows that

\begin{eqnarray*}
\sum_{j \in J''} C \delta^{-1} \rho^3 |z-z_j|^{-2} &\leq& C \delta^{-1} \rho \int_{0}^{2\pi} \int_{(1-\sqrt{2}/2)\rho}^{\operatorname{diam}(K_\delta)} \frac{1}{s^2}s \, ds \, d\theta\\
&\leq& C \delta^{-1} \rho \log \frac{ \operatorname{diam}(K_\delta)}{\rho},
\end{eqnarray*}
which is less than $\epsilon$ provided $\rho$ is sufficiently small. This holds for all $z \in \cup_{j \in J} \partial \mathbb{D}(z_j,r)$, and therefore for all $z \notin \cup_{j \in J} \mathbb{D}(z_j,r)$ by the maximum principle.

Thus, choosing $\rho$ sufficiently small and $N$ sufficiently large, we get
$$\left| \sum_{j \in J} \frac{a_j}{z-z_j} \right| \leq 1+2\epsilon \qquad (z \notin K(Z,R)),$$
where $Z=(z_1,\dots,z_n)$ is the vector of the centers of the squares $Q_j$ with $j \in J$, and $R=(r,\dots,r).$ The result then follows by observing that
$$\gamma(K)=f'(\infty) = (f\phi)'(\infty) = \frac{1}{\pi} \int f(w) \overline{\partial}\phi(w) \, dA(w) = \sum_{j \in J}a_j,$$
where we used Equation (\ref{eqCG}).

\end{proof}

We can now prove Theorem \ref{disks}.

\begin{proof}

Suppose that there exist compact sets $E,F$ with
$$
\gamma(E \cup F) > \gamma(E) + \gamma(F).
$$
Let $0<\epsilon < \gamma(E \cup F) - \gamma(E) - \gamma(F)$. Take $\delta >0$ sufficiently small so that
\begin{equation}
\label{eqqq1}
\gamma(E_{2\delta}) < \gamma(E) + \epsilon/3,
\end{equation}
and
\begin{equation}
\label{eqqq2}
\gamma(F_{2\delta}) < \gamma(F) + \epsilon/3.
\end{equation}
By Lemma \ref{LemMelnikov}, there exist $z_1, z_2, \dots, z_n \in (E \cup F)_{\delta}$ and $0<r<\delta$ such that
$$
\gamma(K(Z,R)) \geq \gamma(E \cup F) - \epsilon/3
$$
and
the disks $\overline{\mathbb{D}}(z_j,r)$ are pairwise disjoint. For each $j \in \{1,2, \dots, n\}$, fix $w_j \in E \cup F$ with $|z_j-w_j|=\operatorname{dist}(z_j,E \cup F) \leq \delta$. Let $A$ be the union of the disks $\overline{\mathbb{D}}(z_j,r)$ with $w_j \in E$, and let $B$ be the union of the disks $\overline{\mathbb{D}}(z_k,r)$ with $w_k \in F \setminus E$.
Then $A \subset E_{2\delta}$ and $B \subset F_{2\delta}$. Finally, we get
\begin{align*}
\gamma(A \cup B) &\geq \gamma(E \cup F) - \epsilon/3\\
&> \gamma(E) + \gamma(F) + \epsilon - \epsilon/3\\
&= \gamma(E) + \epsilon/3 + \gamma(F) + \epsilon/3\\
&> \gamma(A) + \gamma(B),
\end{align*}
where we used equations (\ref{eqqq1}) and (\ref{eqqq2}).

\end{proof}

\subsection{Numerical experiments}

We saw in the previous subsection that in Problem \ref{prob1}, it suffices to consider compact sets $E$ and $F$ that are disjoint finite unions of disjoint closed disks, all with the same radius. In other words, if there is a counterexample to the subadditivity inequality, then there must be a counterexample with disks. This is quite fortunate, since our numerical method for the computation of analytic capacity converges very quickly for such sets, see Example \ref{ex1} and Example \ref{ex2}. This allows us to perform several numerical experiments, in the hope of perhaps finding a counterexample.

To describe these numerical experiments, let $z_1, \dots, z_n, w_1, \dots, w_m$ be distinct complex numbers in $\mathbb{C}$, and let $\delta>0$ be the minimal distance separating these points. For $0<r<\delta/2$, the disks $\overline{\mathbb{D}}(z_j,r)$, $j=1, \dots, n$, $\overline{\mathbb{D}}(w_k,r)$, $k=1,\dots,m$ are all pairwise disjoint. In this case, we set $E_r := \cup_j \overline{\mathbb{D}}(z_j,r)$ and $F_r:= \cup_k \overline{\mathbb{D}}(w_k,r)$.

By Theorem \ref{disks}, the subadditivity of analytic capacity is equivalent to the inequality
$$R(r):= \frac{\gamma(E_r \cup F_r)}{\gamma(E_r)+\gamma(F_r)} \leq 1$$
for all $z_1, \dots, z_n, w_1, \dots, w_m$ and all $0<r<\delta/2$.

Using the algorithm of Subsection \ref{subsec32}, one can compute, for given $z_1,\dots,z_n$ and $w_1,\dots,w_m$, the ratio $R(r)$ for many values of $r$, and then plot the graph of $R(r)$ versus $r$.
\newpage
\begin{example}
\label{ex5}
In this example, the total number of disks is $40$. The compact set $E$ is composed of the $20$ disks with bold boundaries, and $F$ is the union of the remaining disks.

\begin{figure}[!h]
\begin{center}
\includegraphics[width=6cm, height=6cm]{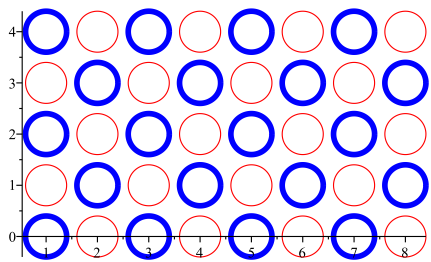}
\caption{The set $E_r \cup F_r$ for Example \ref{ex5}}
\end{center}
\end{figure}

For $500$ values of the radius $r$ equally distributed between $0$ and $0.499$, we computed lower and upper bounds for the ratio $R(r)$.
Figure~\ref{F:40disksgraph} shows the graph of the lower bound versus $r$. The graph for the upper bound is almost identical; the two graphs differ by at most $0.002481$.

\begin{figure}[!h]
\begin{center}
\includegraphics[width=7cm, height=7cm]{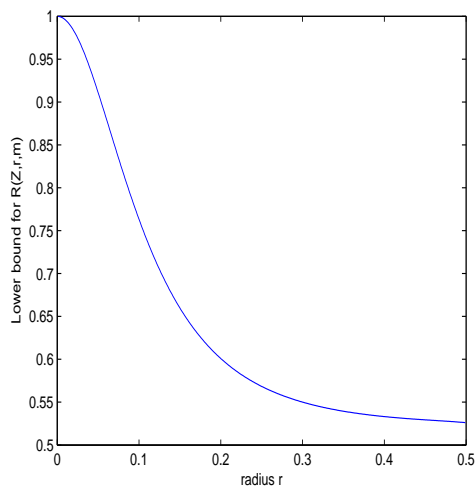}
\caption{Graph of the ratio $R(r)$ for Example \ref{ex5}}
\label{F:40disksgraph}
\end{center}
\end{figure}

\end{example}

\begin{example}
\label{ex6}
In this example, the total number of disks is $18$. The compact set $E$ is composed of the $12$ disks with bold boundaries, and $F$ is the union of the remaining disks.

\begin{figure}[!h]
\begin{center}
\includegraphics[width=7cm, height=7cm]{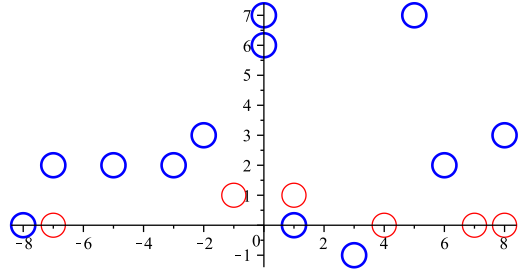}
\caption{The set $E_r \cup F_r$ for Example \ref{ex6}}
\end{center}
\end{figure}

\begin{figure}[!h]
\begin{center}
\includegraphics[width=7cm, height=7cm]{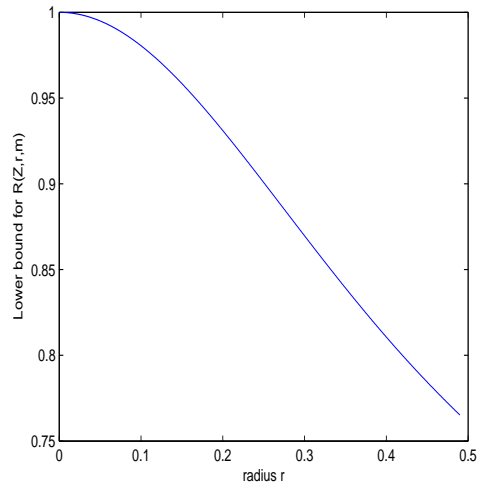}
\caption{Graph of the ratio $R(r)$ for Example \ref{ex6}}
\end{center}
\end{figure}

\end{example}

A few remarks are in order. First, all of our numerical experiments seem to suggest that $R(r) \leq 1$ for all $r$, i.e. that analytic capacity is indeed subadditive. More surprisingly though, all computations seem to suggest that the following might be true.

\begin{conjecture}[Younsi--Ransford \cite{YOU}]
\label{conjdecreasing}
For all $z_1, \dots, z_n, w_1, \dots, w_m$, the function $R(r)$ is decreasing in $0<r<\delta/2$.
\end{conjecture}

A proof of Conjecture \ref{conjdecreasing} would imply that analytic capacity is subadditive. Indeed, this follows from the following asymptotic expression for $R(r)$ as $r \to 0^{+}$, which shows that for fixed centers of the disks, the ratio $R(r)$ is less than 1, for all $r$ sufficiently small.

\begin{theorem}[Younsi--Ransford \cite{YOU}]
\label{asymp}
Let $z_1, \dots, z_n, w_1, \dots, w_m \in \mathbb{C}$. Then
$$R(r)=1-Cr^2+O(r^3) \qquad (r \to 0^{+}),$$
where $C$ is a positive constant depending only on $m$, $n$ and the points $z_1, \dots, z_n$ and $w_1, \dots, w_m$.

\end{theorem}

The proof also follows from results of Melnikov \cite{MEL}. The constant $C$ and the error term in Theorem \ref{asymp} can be made more or less explicit, but it seems difficult to extract any useful information.

Lastly, we mention that Conjecture \ref{conjdecreasing} was proved in the simplest case $n=m=1$ (see \cite[Theorem 8.2]{YOU}), based on Murai's formula from \cite{MUR2} for the analytic capacity of the union of two disjoint closed disks.

\section{Rational Ahlfors functions}
\label{sec5}

In this section, we discuss another problem related to the computation of analytic capacity.

\subsection{Definition and properties}
\label{subsec51}

The starting point of the study of rational Ahlfors functions is the following result of Jeong and Taniguchi \cite{JEO}, whose uniqueness part was proved in \cite{FOR}. Recall that a domain $\Omega$ in the Riemann Sphere $\RiemannSphere$ is a non-degenerate $n$-connected domain if $\RiemannSphere \setminus \Omega$ has exactly $n$ connected components, each of them containing more than one point. By a repeated application of the Riemann mapping theorem, every such domain $\Omega$ is conformally equivalent to a domain bounded by analytic curves. We may therefore assume without loss of generality that $\partial \Omega$ is analytic.

\begin{theorem}[Jeong--Taniguchi \cite{JEO}, Fortier Bourque--Younsi \cite{FOR}]
\label{JT}
Let $\Omega$ be a non-degenerate $n$-connected domain containing $\infty$, and let $f$ be the Ahlfors function on $\Omega$. Then there exist a rational map $R$ of degree $n$ and a conformal map $g:\Omega \to R^{-1}(\mathbb{D})$ such that $f=R \circ g$. Moreover, if $f=Q \circ h$ for another rational map $Q$ of degree $n$ and a conformal map $h: \Omega \to Q^{-1}(\mathbb{D})$, then there is a M\"{o}bius transformation $T$ such that $Q=R \circ T^{-1}$ and $h=T \circ g$.
\end{theorem}

In particular, this shows that every non-degenerate $n$-connected domain is conformally equivalent to a domain of the form $R^{-1}(\mathbb{D})$ for some rational map $R$.

The proof relies on the construction of a Riemann surface obtained by welding $n$ copies of the unit disk $\mathbb{D}$ to $\Omega$ using the map $f$. The Riemann surface $X$ thereby obtained is topologically a sphere, thus there is a conformal map $g: X \to \RiemannSphere$. One can then easily extend $f$ to $X$ so that the composition $f \circ g^{-1} : \RiemannSphere \to \RiemannSphere$ is holomorphic and hence a rational map $R$ of the degree $n$, the same degree as the proper map $f$. This proves the existence statement, and the uniqueness part follows from a standard removability argument.

The uniqueness part of Theorem \ref{JT} implies in particular that the pair $(R,g)$ is unique if we require the conformal map $g:\Omega \to R^{-1}(\mathbb{D})$ to be normalized by $g(z)=z+b_1/z + b_2/z^2 + \dots$ near $\infty$. Moreover, in this case, the conformal invariance of Ahlfors functions (see Subsection \ref{subsec22}) shows that $R$ is the Ahlfors function on its restricted domain $R^{-1}(\mathbb{D})$. Such rational maps are called rational Ahlfors functions.
More precisely, a rational map $R$ is a \textit{rational Ahlfors function} if

\begin{itemize}
\item $R(\infty)=0$
\item $R^{-1}(\mathbb{D})$ is an $n$-connected domain
\item $R$ is the Ahlfors function on $R^{-1}(\mathbb{D})$.
\end{itemize}

The second condition is equivalent to requiring that all the critical values of $R$ belong to $\mathbb{D}$. In particular, a rational Ahlfors function must have only simple poles, since otherwise it would have $\infty$ as a critical value, and can therefore be written as
$$R(z)=\sum_{j=1}^n \frac{a_j}{z-p_j}$$
for some $a_1,\dots,a_n \in \mathbb{C} \setminus \{0\}$ and distinct $p_1,\dots,p_n \in \mathbb{C}$. With this representation, we have $R'(\infty)=\sum_{j=1}^n a_j$, so that $R$ is a rational Ahlfors function if and only if $\gamma(E)=\sum_{j=1}^n a_j$, where $E:=\RiemannSphere \setminus R^{-1}(\mathbb{D})=R^{-1}(\RiemannSphere \setminus \mathbb{D})$.

In \cite{JEO}, Jeong and Taniguchi raised the problem of finding all rational Ahlfors functions.

\begin{problem}
\label{prob2}
Determine which residues $a_1,\dots,a_n$ and poles $p_1,\dots,p_n$ correspond to rational Ahlfors functions.
\end{problem}

In view of Theorem \ref{JT}, a solution to Problem \ref{prob2} would yield a complete understanding of all Ahlfors functions on finitely connected domains, up to conformal equivalence.

We mention that the problem is trivial for $n=1$ : rational Ahlfors functions of degree $1$ are precisely the rational maps of the form $z \mapsto a_1/(z-p_1)$ for $a_1>0$ and $p_1 \in \mathbb{C}$. As for higher degree, the first examples were given in \cite{FOR}.

\subsection{Rational Ahlfors functions with reflection symmetry}
\label{subsec52}

Let $R$ be a rational map with $R(\infty)=0$, and suppose that $\Omega:=R^{-1}(\mathbb{D})$ is an $n$-connected domain. In this case, the domain $\Omega$ is bounded by $n$ disjoint analytic curves $\Gamma_1, \dots, \Gamma_n$. Let $f$ be the Ahlfors function on $\Omega$. By definition, proving that $R$ is a rational Ahlfors function amounts to showing that $R=f$.

Now, both $R$ and $f$ are degree $n$ proper analytic maps of $\Omega$ onto the unit disk $\mathbb{D}$. In particular, they extend analytically across $\partial \Omega$ and map each boundary curve $\Gamma_j$ homeomorphically onto $\partial \mathbb{D}$. In general, however, there are many such proper maps, in view of the following result of Bieberbach and Grunsky, whose proof can be found in \cite[Theorem 2.2]{BELL} or \cite[Theorem 3]{KHA}.

\begin{theorem}[Bieberbach \cite{BIEB}, Grunsky \cite{GRU}]
\label{bieber}
Let $\Omega$ be a non-degenerate $n$-connected domain containing $\infty$ and bounded by $n$ disjoint analytic Jordan curves $\Gamma_1, \dots, \Gamma_n$. For each $j$, let $\alpha_j$ be any point in $\Gamma_j$. Then there exists a unique proper analytic map $g: \Omega \to \mathbb{D}$ of degree $n$ satisfying $g(\infty)=0$ whose extension to $\partial \Omega$ maps each $\alpha_j$ to the point $1$.
\end{theorem}

In particular, one way to prove that $R=f$ is to show that these two maps send the same points to $1$. This will be the case if, for instance, $R$ has only real poles and positive residues.

\begin{theorem}[Fortier Bourque--Younsi \cite{FOR}]
\label{theosym1}
Let
$$R(z):=\sum_{j=1}^n \frac{a_j}{z-p_j},$$
where the poles $p_1,\dots,p_n$ are distinct and real, and the residues $a_1,\dots,a_n$ are positive. If $\Omega:=R^{-1}(\mathbb{D})$ is an $n$-connected domain, then $R$ is a rational Ahlfors function.
\end{theorem}
The idea of the proof is to observe that in this case, the domain $\Omega$ is symmetric with respect to the real axis, intersecting it at $2n$ points $\alpha_1<\beta_1<\dots<\alpha_n<\beta_n$ such that $R(\alpha_j)=-1$ and $R(\beta_j)=1$ for each $j$. Now, the domain $\Omega$ can be mapped conformally onto the complement $\Omega'$ in $\RiemannSphere$ of $n$ disjoint closed intervals in the real line, by mapping $\Omega \cap \mathbb{H}$ onto the upper half-plane $\mathbb{H}$, and then using the Schwarz reflection principle. Now, using the formula for the Ahlfors function on $\Omega'$ and conformal invariance (see Subsection \ref{subsec22}), one can show that $f(\beta_j)=1$ for each $j$, from which it follows that $R=f$ by the uniqueness part of Theorem \ref{bieber}.

\subsection{Rational Ahlfors functions with rotational symmetry}
\label{subsec53}

Again, let $R$ be a rational map with $R(\infty)=0$, and assume that $\Omega=R^{-1}(\mathbb{D})$ is an $n$-connected domain. Let $f : \Omega \to \mathbb{D}$ be the Ahlfors function. We saw in Subsection \ref{subsec52} that in order to prove that $R=f$, it suffices to show that they both map the same points to $1$. Another way to prove that these maps are equal is to show that they have the same zeros, because then one can apply the maximum principle to both quotients $R/f$ and $f/R$, whose absolute values are equal to $1$ everywhere on $\partial \Omega$, to deduce that $R/f$ is a unimodular constant, which has to be one provided $R'(\infty)>0$. This observation was used in \cite{FOR} to prove the following result.

\begin{theorem}[Fortier Bourque--Younsi \cite{FOR}]
\label{theosym2}
Let $n\geq 2$, $0<a < n (n-1)^{(1-n)/n}$, and $R(z)=a z^{n-1} / (z^n - 1)$. Then $R$ is a rational Ahlfors function.
\end{theorem}

The condition on $a$ guarantees that $\Omega$ is $n$-connected. Now, if $\omega:=e^{2\pi i/n}$, then $\omega \Omega = \Omega$ and thus
$$\omega f(\omega z) = f(z) \qquad (z \in \Omega),$$
by uniqueness of the Ahlfors function. This implies that $f$ vanishes only at $0$ and $\infty$, the same zeros as $R$. Indeed, if $f$ had a zero at $z_0 \neq 0,\infty$, then $f$ would vanish at $z_0, \omega z_0,\dots, \omega^{n-1}z_0$ and $\infty$, a total of $n+1$ distinct points, contradicting the fact that $f$ has degree $n$. It follows that $R=f$.

\subsection{Positivity of residues and numerical examples}
\label{subsec54}

As mentioned at the end of Subsection \ref{subsec51}, the rational Ahlfors functions of degree one are of the form $R(z)=a/(z-p)$, where $a>0$ and $p \in \mathbb{C}$. In degree two, the positivity of residues is also a necessary and sufficient condition for a rational map to be a rational Ahlfors function, as proved in \cite{FOR}.

\begin{theorem}[Fortier Bourque--Younsi \cite{FOR}]
A rational map of degree two is a rational Ahlfors function if and only if it can be written in the form
$$R(z)= \frac{a_1}{z-p_1}+\frac{a_2}{z-p_2},$$
for distinct $p_1, p_2 \in \mathbb{C}$ and positive $a_1,a_2$ satisfying $a_1+a_2 < |p_1-p_2|$.
\end{theorem}

Now, note that all the examples of rational Ahlfors functions that we obtained so far have positive residues. It thus seems natural to expect the positivity of residues to be a sufficient condition, in any degree. Unfortunately, as observed in \cite{FOR}, this fails even in degree $3$. Before presenting the counterexample, we first explain how to numerically check whether a given rational map is a rational Ahlfors function.

Let
$$R(z):=\sum_{j=1}^n \frac{a_j}{z-p_j},$$
and suppose that $R^{-1}(\mathbb{D})$ is an $n$-connected domain. Let $E$ be the compact set
$$E:=R^{-1}(\RiemannSphere \setminus \mathbb{D}) = \{z \in \mathbb{C} : |R(z)| \geq 1\}.$$
Recall that $R$ is a rational Ahlfors function if and only if $\gamma(E) = \sum_{j=1}^\infty a_j$. We can therefore check whether $R$ is a rational Ahlfors function by computing $\gamma(E)$ using the numerical method of Subsection \ref{subsec32}. However, the method involves integrals over the boundary of $E$ with respect to arclength, and thus the first step is to obtain a parametrization of $\partial E = \{z \in \mathbb{C} : |R(z)|=1\}$. One can do this by solving for $z$ in the equation $R(z)=e^{it}$, which is easily done provided the degree $n$ of $R$ is small, say $n \leq 3$.

We now present some numerical examples to illustrate the method. All the numerical work was done with \textsc{matlab}, the integrals being computed with a precision of $10^{-9}$.

\newpage
\begin{example}
\label{ex7}
Let
$$R(z) = \frac{0.2}{z+2} + \frac{0.1}{z} + \frac{0.4}{z-5}.$$

\begin{figure}[h!t!b]
\begin{center}
\includegraphics[width=6cm]{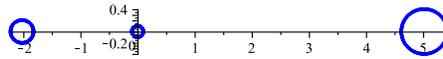}
\caption{The boundary of $E=\{\,z \in \mathbb{C}: |R(z)|\geq 1 \,\}$ for Example \ref{ex7}}
\end{center}
\end{figure}

As described in Subsection \ref{subsec32}, one can compute $\gamma(E)$ using linear combinations of the functions $(z+2)^{-j}$, $z^{-j}$ and $(z-5)^{-j}$, $1 \leq j \leq k$, for various values of $k$.

\begin{table}[!hbp]
\begin{center}
\caption{Lower and upper bounds for $\gamma(E)$ for Example \ref{ex7}}
\label{table3}
\begin{tabular}{|c|l|l|}
\hline
$k$ & Lower bound for $\gamma(E)$ & Upper bound for $\gamma(E)$ \\
\hline
$1$ &   0.696735209508754 &     0.700011861859377\\
$2$ &    0.699988138057939 &       0.700000163885012 \\
$3$ &     0.699999835775098 &     0.700000002518033 \\
\hline
\end{tabular}
\end{center}
\end{table}

By Theorem \ref{theosym1}, the rational map $R$ is the Ahlfors function for the compact set $E$, and we have
$$\gamma(K) = R'(\infty) = 0.2+0.1+0.4=0.7.$$
Our numerical results therefore agree with the predicted value.

\end{example}
\newpage
\begin{example}
\label{ex8}
Let
$$R(z) = \frac{z^2}{z^3-1}.$$

\begin{figure}[h!t!b]
\begin{center}
\includegraphics[width=6cm, height=6cm]{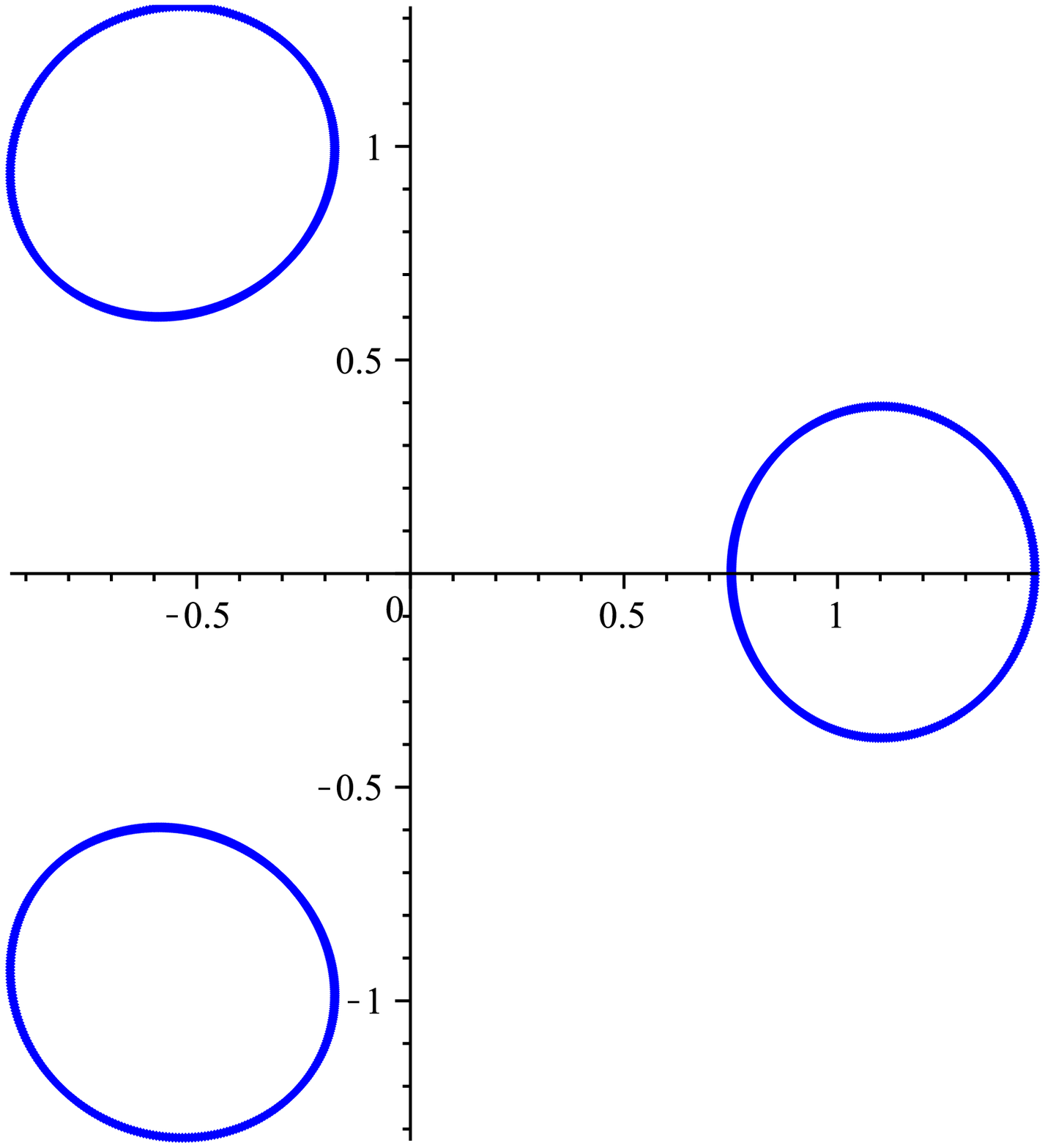}
\caption{The boundary of $E=\{\,z \in \mathbb{C}: |R(z)|\geq 1 \,\}$ for Example \ref{ex8}}
\end{center}
\end{figure}

\begin{table}[!hbp]
\begin{center}
\caption{Lower and upper bounds for $\gamma(E)$ for Example \ref{ex8}}
\label{table4}
\begin{tabular}{|c|l|l|}
\hline
$k$ & Lower bound for $\gamma(E)$ & Upper bound for $\gamma(E)$ \\
\hline
$1$ &   0.897012961211562 &      1.003766600572323\\
$2$ &    0.996247533470256&       1.000449247199905 \\
$3$ &     0.999550954532515 &     1.000227970885994 \\
$4$ &     0.999772081072887 &      1.000015305500631 \\
$5$ &     0.999984694733624 &     1.000004234543914 \\
$6$ &     0.999995765474017 &      1.000002049275081 \\
\hline
\end{tabular}
\end{center}
\end{table}

By Theorem \ref{theosym2}, $R$ is the Ahlfors function for the compact set $E$, and thus we have
$$\gamma(E) = R'(\infty) = 1.$$
\end{example}
\newpage
\begin{example}
\label{ex9}
The following example answers a question raised in \cite{FOR} asking whether Theorem \ref{theosym1} holds if the rational map $R$ is allowed to have conjugate pairs of poles.

\begin{question}[Fortier Bourque--Younsi \cite{FOR}]
\label{q1}
Suppose that a rational map $R$ satisfies $R(\overline z) = \overline{R(z)}$ and has all positive residues. If $R^{-1}(\mathbb{D})$ is an $n$-connected, must $R$ be a rational Ahlfors function?
\end{question}

The following numerical counterexample shows that the answer to Question \ref{q1} is in fact negative. Consider

$$R(z)=\frac{0.6}{z}+\frac{1.2}{z-(1+4i)}+\frac{1.2}{z-(1-4i)}.$$

\begin{figure}[!h]
\begin{center}
\includegraphics[width=8cm, height=8cm]{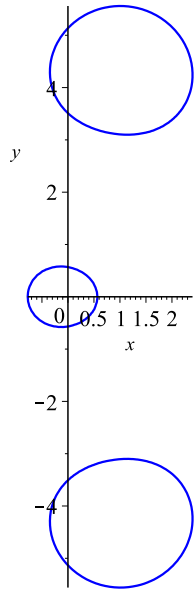}
\caption{The boundary of $E=\{\,z \in \mathbb{C}: |R(z)|\geq 1 \,\}$ for Example \ref{ex9}}
\end{center}
\end{figure}

\begin{table}[!hbp]
\begin{center}
\caption{Lower and upper bounds for $\gamma(E)$ for Example \ref{ex9}}
\label{table4}
\begin{tabular}{|c|l|l|}
\hline
$k$ & Lower bound for $\gamma(E)$ & Upper bound for $\gamma(E)$ \\
\hline
$1$ &   2.791098712682427 &      3.011700281953340\\
$2$ &    2.990065509632496&       3.003544766713681 \\
$3$ &     2.998352995858360 &     3.001339199421414 \\
$4$ &     3.000606967107599 &       3.001054476070951 \\
$5$ &     3.000898172927223 &     3.000989499718515 \\
$6$ &     3.000964769941962 &      3.000979603389716 \\
$7$ & 3.000974855197433 & 3.000977619119024\\
\hline
\end{tabular}
\end{center}
\end{table}
We see that $R'(\infty) = 3 < 3.000974855197433 \leq \gamma(E)$, so that $R$ is not a Rational Ahlfors function.

\end{example}

The above example also shows that the positivity of residues is not sufficient for a rational map of degree $3$ to be Ahlfors. This can be also be extended to any degree $n \geq 3$, replacing $R$ by
$$R_\epsilon(z):= R(z) + \sum_{j=4}^n \frac{\epsilon}{z-b_j},$$
where $b_4,\dots, b_n$ are distinct points in $R^{-1}(\mathbb{D})$ and $\epsilon>0$ is small.

\begin{theorem}[Fortier Bourque--Younsi \cite{FOR}]
For every $n \geq 3$, there exists a rational map $R$ such that $R(\infty)=0$, $R^{-1}(\mathbb{D})$ is an $n$-connected domain, and $R$ has only positive residues, but $R$ is not a rational Ahlfors function.
\end{theorem}

Finally, using Koebe's continuity method based on Brouwer's Invariance of Domain theorem, it is possible to show that the residues of a rational Ahlfors functions are not necessarily positive.

\begin{theorem}[Fortier Bourque--Younsi \cite{FOR}]
\label{notnec}
For every $n \geq 3$, there exists a rational Ahlfors function of degree $n$ whose residues are not all positive.
\end{theorem}
Combining these results show that the positivity of residues is neither sufficient for necessary for a rational map to be Ahlfors, in any degree $n \geq 3$.

The proof of Theorem \ref{notnec}, however, is not constructive, and it would be very interesting to find explicit examples.

\begin{problem}
\label{prob3}
Find an explicit example of a rational Ahlfors function whose residues are not all positive.
\end{problem}

\section{The Cauchy capacity}
\label{sec6}

This section deals with the relationship between analytic capacity and the Cauchy capacity, another similar extremal problem. Before proceeding further, let us first give some motivation with a brief overview of Tolsa's solution of Painlev\'e's problem.

\subsection{The capacity $\gamma_+$, Painlev\'es problem and the semi-additivity of analytic capacity}
\label{subsec61}

We only give a very brief introduction to Painlev\'es problem and related results, since a detailed discussion is beyond the scope of this survey article. The interested reader may consult \cite{DUD}, \cite{TOL} or \cite{TOL2} for more information.

Recall from the introduction that Painlev\'e's problem asks for a geometric characterization of the compact sets that are removable for bounded analytic functions. We say that a compact set $E \subset \mathbb{C}$ is \textit{removable} (for bounded analytic functions) if every bounded analytic function on $\mathbb{C} \setminus E$ is constant. As mentioned in Subsection \ref{subsec21}, removable sets coincide precisely with the sets of zero analytic capacity.

Painlev\'e was the first one to observe that there is a close relationship between removability and Hausdorff measure and dimension. More precisely, he proved that compact sets of finite one-dimensional Hausdorff measure are removable. In particular, sets of dimension less than one are removable. On the other hand, any compact set $E$ with Hausdorff dimension bigger than one is not removable. This follows from Frostman's lemma, which gives the existence of a non-trivial Radon measure $\mu$ supported on $E$ whose \textit{Cauchy transform}

$$\mathcal{C}\mu (z):=\int \frac{1}{\xi - z} d\mu(\xi)$$
is continuous on $\RiemannSphere$. In particular, $\mathcal{C}\mu$ is bounded, and, since it is analytic on $\mathbb{C} \setminus E$ and non-constant ($\mathcal{C}\mu(\infty)=0$ but $\mathcal{C}\mu'(\infty) = -\mu(E) \neq 0$), we get that $E$ is not removable.

This proof illustrates the important role played by the Cauchy transform in the study of removable sets, mainly because it is an easy way to construct functions analytic outside a given compact set.

In view of the above remarks, Painlev\'e's problem reduces to the case of dimension exactly equal to one. In dimension one though, the problem quickly appeared to be extremely difficult, and it took a very long time until progress was made. One of the major advances was the proof of the so-called Vitushkin's conjecture by David \cite{DAV} in 1998.

\begin{theorem}[Vitushkin's conjecture]
Let $E$ be a compact set with finite one-dimensional Hausdorff measure. Then $E$ is removable if and only if $E$ is purely unrectifiable, i.e. it intersects every rectifiable curve in a set of zero one-dimensional Hausdorff measure.
\end{theorem}

The forward implication was previously known as Denjoy's conjecture and in fact follows from the results of Calder\'on on the $L^2$-boundedness of the Cauchy transform operator. We also mention that the converse implication is false without the assumption that $E$ has finite length, see \cite{JOM}

In order to obtain a characterization of removable sets in the case of dimension one but infinite length, Tolsa \cite{TOL3} proved that analytic capacity is comparable to a quantity which is easier to comprehend as it is more suitable to real analysis tools. More precisely, define the capacity $\gamma_{+}$ of a compact set $E$ by
$$\gamma_{+}(E):=\sup \{\mu(E) : \operatorname{supp}{\mu} \subset E, |\mathcal{C}\mu| \leq 1 \, \, \mbox{on}\, \, \RiemannSphere \setminus E\},$$
where $\mu$ is a positive Radon measure supported on $E$. Clearly, we have $\gamma_{+}(E) \leq \gamma(E)$. Note also that the supremum in the definition of $\gamma_+$ is always attained by some measure, by a standard $\operatorname{weak}^*$ convergence argument.

\begin{theorem}[Tolsa \cite{TOL3}]
\label{TheoTolsa}
There is a universal constant $C$ such that
$$\gamma(E) \leq C \gamma_{+}(E).$$
\end{theorem}

This remarkable result has several important consequences. For instance, it gives a complete solution to Painlev\'e's problem for arbitrary compact sets, involving the notion of curvature of a measure introduced by Melnikov \cite{MEL}.

\begin{theorem}[Tolsa \cite{TOL3}]
A compact set $E \subset \mathbb{C}$ is not removable if and only if it supports a nontrivial positive Radon measure with linear growth and finite curvature.
\end{theorem}

We say that a positive Radon measure $\mu$ has \textit{linear growth} if there exists a constant $C$ such that $\mu(\mathbb{D}(z,r)) \leq Cr$ for all $z \in \mathbb{C}$ and all $r >0$. The \textit{curvature} $c(\mu)$ of $\mu$ is defined by
$$c(\mu)^2:=\int \int \int \frac{1}{R(x,y,z)^2} d\mu(x) d\mu(y) d\mu(z),$$
where $R(x,y,z)$ is the radius of the circle passing through $x,y,z$.

Another importance consequence of Theorem \ref{TheoTolsa} is the semi-additivity of analytic capacity, as mentioned in Section \ref{sec4}.

\subsection{The Cauchy capacity}
\label{subsec62}

We can now define the Cauchy capacity. Given a compact set $E \subset \mathbb{C}$, the \textit{Cauchy capacity} of $E$, noted by $\gamma_c(E)$, is defined by
$$\gamma_c(E):= \sup \{|\mu(E)| : \operatorname{supp}{\mu} \subset E, |\mathcal{C}\mu| \leq 1 \, \, \mbox{on}\, \, \mathbb{C}_\infty \setminus E\},$$
where $\mu$ is a complex Borel measure supported on $E$. In other words, the Cauchy capacity is defined in the same way as $\gamma_+$, except that complex measures are allowed. Note that
$$\gamma_+(E) \leq \gamma_c(E) \leq \gamma(E) \leq C \gamma_{+}(E) \leq C \gamma_c(E)$$
for any compact set $E$, where $C$ is the constant of Theorem \ref{TheoTolsa}. In particular, the capacities $\gamma, \gamma_+, \gamma_c$ vanish simultaneously, which is already a deep result.

As far as we know, the following question was raised by Murai \cite{MUR}. See also \cite{HAV2}, \cite{HAV3} and \cite[Section 5]{TOL2}.

\begin{problem}
\label{prob4}
Is analytic capacity actually equal to the Cauchy capacity? In other words, is it true that
\begin{equation}
\label{eqprob4}
\gamma(E)=\gamma_c(E)
\end{equation}
for all compact sets $E \subset \mathbb{C}$?

\end{problem}

In other words, Problem \ref{prob4} asks whether the supremum in the definition of analytic capacity remains unchanged if we only consider bounded analytic functions which are Cauchy transforms of complex measures supported on the set. In particular, Equation (\ref{eqprob4}) holds for $E$ if every bounded analytic function on $\mathbb{C} \setminus E$ vanishing at $\infty$ is the Cauchy transform of a complex measure supported on $E$. This is the case if, for instance, $E$ has finite one-dimensional Hausdorff measure, or, more generally, if it has \textit{finite Painlev\'e length}, meaning that there is a number $l$ such that every open set $U$ containing $E$ contains a cycle $\Gamma$ surrounding $E$ that consists of finitely many disjoint analytic Jordan curves and has length less than $l$. The fact that $\gamma = \gamma_c$ for such sets is easily derived from Cauchy's integral formula and a $\operatorname{weak}^*$ convergence argument.

\begin{proposition}
\label{PropPain}
If $E$ has finite Painlev\'e length, then $\gamma(E)=\gamma_c(E)$.
\end{proposition}

See \cite{YOU3} for a generalization to compact sets of $\sigma$-finite Painlev\'e length, in a sense.

Note that every compact set in the plane can be obtained as a decreasing sequence of compact sets with finite Painlev\'e length. In particular, by outer regularity of analytic capacity, a positive answer to Problem \ref{prob4} would follow if one could prove that the Cauchy capacity is also outer regular.

\begin{problem}
\label{prob5}
Is it true that if $E_n \downarrow E$, then $\gamma_c(E_n) \downarrow \gamma_c(E)$?
\end{problem}

\subsection{Is the Ahlfors function a Cauchy transform?}

In Proposition \ref{PropPain}, not only the Ahlfors function but every bounded analytic function on $\mathbb{C} \setminus E$ vanishing at $\infty$ is the Cauchy transform of a complex measure supported on $E$. From the point of view of Problem \ref{prob4}, a more interesting question is whether the Ahlfors function can always be expressed as the Cauchy transform of a complex measure supported on the set. This was, however, answered in the negative by Samokhin.

\begin{theorem}[Samokhin \cite{SAM}]
\label{TheoSam}
There exists a connected compact set $F$ with connected complement such that the Ahlfors function for $F$ is not the Cauchy transform of any complex Borel measure supported on $F$.
\end{theorem}

In particular, this implies that $\gamma_+(F)<\gamma(F)$. Indeed, suppose that $\gamma_+(F)=\gamma(F)$, and let $\mu$ be a positive Borel measure supported on $F$ with $|\mathcal{C}\mu| \leq 1$ on $\mathbb{C} \setminus F$ and $\mu(F)=\gamma_+(F)$. Then $g:=\mathcal{C}(-\mu)$ is analytic on $\mathbb{C} \setminus F$ and satisfies $|g| \leq 1$ on $\mathbb{C} \setminus F$. But $g'(\infty)=\mu(F)=\gamma_+(F)=\gamma(F)$, so that $g$ is the Ahlfors function for $F$, by uniqueness, contradicting Theorem \ref{TheoSam}.

This argument relies on the fact that the supremum in the definition of $\gamma_+$ is always attained by some measure, which follows from a standard $\operatorname{weak}^*$ convergence argument. It is not clear a priori whether this remains true for the Cauchy capacity $\gamma_c$, since in this case one has to deal with complex measures.

In fact, in \cite{YOU3}, we constructed a compact set $E$ for which there is no complex Borel measure $\mu$ supported on $E$ such that
$$|\mathcal{C}\mu(z)| \leq 1 \qquad (z \in \mathbb{C} \setminus E)$$
and $\mu(E)=\gamma_c(E)$. This follows from the following result, by the same argument as above.

\begin{theorem}[Younsi \cite{YOU3}]
\label{TheoCauchyCap}
There exists a connected compact set $E$ with connected complement such that $\gamma(E)=\gamma_c(E)$, but the Ahlfors function for $E$ is not the Cauchy transform of any complex Borel measure supported on $E$.
\end{theorem}

The construction is a bit simpler than the one in \cite{SAM}, which makes it easier to show that the analytic capacity and the Cauchy capacity of the set are equal. The set $E$ is the union of the nonrectifiable curve $\Gamma:=\{x + ix\sin{(1/x)} : x \in (0,1/\pi]\}$ and the line segment $[-i,i]$.

\begin{figure}[h!t!b]
\begin{center}
\includegraphics[width=6cm, height=6cm]{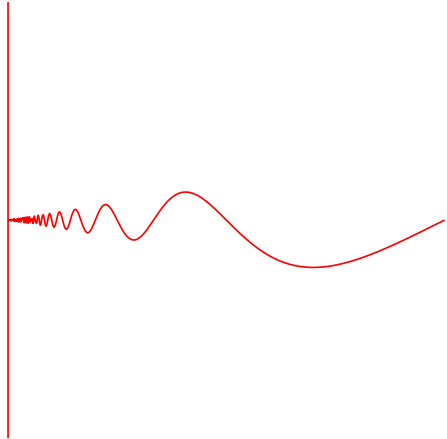}
\caption{The compact set $E$.}
\end{center}
\end{figure}

The proof that the Ahlfors function $f$ for $E$ is not a Cauchy transform relies on the fact that since $E$ is connected, the Ahlfors function is a conformal map of $\RiemannSphere \setminus E$ onto $\mathbb{D}$ (see Subsection \ref{subsec22}). More precisely, if $f=\mathcal{C}\mu$ for some $\mu$ supported on $E$, then one can use Cauchy's formula to relate the measure $\mu$ on $\Gamma \setminus \{0\}$ with the boundary values of $f$. Combining well-known results on boundary correspondence under conformal maps with the fact that the curve $\Gamma$ has infinite length then makes it possible to show that the total variation of $\mu$ must be infinite, a contradiction.

As for the proof that $\gamma(E)=\gamma_c(E)$, it essentially follows from a convergence result for analytic capacity. For $k \in \mathbb{N}$, let $E_k$ be the union of the segment $[-i,i]$ and the portion of $\Gamma$ with $1/(\pi+k)<x \leq 1/\pi$. Since every $E_k$ has finite length, we know that $\gamma(E_k)=\gamma_c(E_k) \leq \gamma_c(E)$ for all $k$, the last inequality because $E_k \subset E$. It thus only remains to prove that $\gamma(E_k) \uparrow \gamma(E)$ as $k \to \infty$. Using conformal invariance, this can be reduced to proving the following lemma.

\begin{lemma}
\label{LemConv}
For $k \in \mathbb{N}$, let $F_k$ be a union of two disjoint closed disks. Suppose that $F_k \to F$, where $F$ is the union of two closed disks intersecting at exactly one point. Then $\gamma(F_k) \to \gamma(F)$.
\end{lemma}

Lemma \ref{LemConv} follows from conformal invariance and the observation that the conclusion holds if we replace disks by intervals contained in the same horizontal line, by Pommerenke's theorem (see Subsection \ref{subsec22}). This completes the sketch of the proof of Theorem \ref{TheoCauchyCap}.

\subsection{Convergence results for analytic capacity}

The proof of Theorem \ref{TheoCauchyCap} shows how even simple questions related to convergence of analytic capacity can be difficult. For example, we do not know whether Lemma \ref{LemConv} remains true if the number of disks is bigger than two.

\begin{problem}
\label{prob6}
Does the conclusion of Lemma \ref{LemConv} still hold if each $F_k$ is instead assumed to be a union of $n$ disks, $n \geq 3$, and one pair of disks intersect at one point in the limit?
\end{problem}

It would be interesting to study this problem numerically using the method of Section \ref{sec3} for the computation of analytic capacity.

Another open problem is the inner-regularity of analytic capacity.

\begin{problem}
\label{prob7}
Suppose that $E_n \uparrow E$, where $E$ and the $E_n$'s are compact. Does $\gamma(E_n) \uparrow \gamma(E)$?
\end{problem}

This should hold if analytic capacity is truly a capacity, in the sense of Choquet. It is rather unfortunate that this remains unknown!

Finally, we end this section by mentioning that it would be interesting to develop a numerical method for the computation of the capacities $\gamma_+$ or $\gamma_c$, similar to the one from \cite{YOU} for analytic capacity. A first step would be to settle the following question which, as far as we know, is still unanswered.

\begin{problem}
\label{prob8}
Is the capacity $\gamma_+$ always attained by a \textit{unique} measure?
\end{problem}


\bibliographystyle{amsplain}

\end{document}